\documentclass[12pt,letterpaper]{amsart}
\usepackage{graphicx}
\usepackage{amsmath, amssymb, amsthm, mathrsfs,verbatim,amscd}
\usepackage[curve]{xypic}
\usepackage[left=1in,top=1in,right=1in]{geometry}
\usepackage{setspace}
\newtheorem{lemma}{Lemma}
\newtheorem{theorem}{Theorem}
\newtheorem{proposition}{Proposition}
\newtheorem{definition}{Definition}
\newtheorem{corollary}{Corollary}
 
\newtheorem{condition}{Condition}
\newtheorem{remark}{Remark}

\newcommand{\R}{\ensuremath{\mathbb{R}} }

\newcommand{\N}{\ensuremath{\mathbb{N}} }

\newcommand{\Z}{\ensuremath{\mathbb{Z}} }

\newcommand{\F}{\ensuremath{\mathcal{F}} }

\newcommand{\T}{\mathcal{T}}

\renewcommand{\S}{\mathcal{S}}
\newcommand{\Lim}{\displaystyle \lim_}
\newcommand{\cf}{\mathbf{1}}
\renewcommand{\P}{\mathbb{P}}
\newcommand{\E}{\mathbb{E}}

\renewcommand{\deg}{\textrm{deg}}

\newcommand{\calP}{\mathcal{P}}
\newcommand{\GW}{\mathrm{GW}}

\newcommand{\gw}{\textrm{GW}}
\newcommand{\bP}{\mathbf{P}}
\newcommand{\deq}{\mathrel{\mathop :}=}
\newcommand{\var}{\textrm{Var}}
\newcommand{\Ht}{\textrm{ht}}
\newcommand{\M}{\mathcal{M}_w}
\newcommand{\frag}{\mathrm{Frag}}

\usepackage{calc}
\usepackage{tikz}
\usetikzlibrary{decorations.markings}
\tikzstyle{vertex}=[circle, draw, inner sep=0pt, minimum size=6pt]
\tikzstyle{Vertex}=[circle, draw, inner sep=0pt, minimum size=14pt]
\newcommand{\vertex}{\node[vertex]}
\newcommand{\Vertex}{\node[Vertex]}

\newcommand{\rspine}{\textrm{rspine}}

\begin{document}

\title[Markov branching and Galton-Watson trees]{Scaling limits of Markov branching trees and Galton-Watson trees conditioned on the number of vertices with out-degree in a given set}

\author{Douglas Rizzolo}
\address{Department of Mathematics \\ University of Washington, WA 98105}
\email{drizzolo@math.washington.edu}
\thanks{This material is based upon work supported in part by the National Science Foundation Graduate 
Research Fellowship under Grant No. DGE 1106400 and in part by NSF DMS-1204840.}

\begin{abstract}
We obtain scaling limits for Markov branching trees whose size is specified by the number of nodes whose out-degree lies in a given set.  This extends recent results of Haas and Miermont in \cite{HaMi10}, which considered the case when the size of a tree is either its number of leaves or its number of vertices.  We use our result to prove that the scaling limit of finite variance Galton-Watson trees conditioned on the number of nodes whose out-degree lies in a given set is the Brownian continuum random tree.  The key to applying our result for Markov branching trees to conditioned Galton-Watson trees is a generalization of the classical Otter-Dwass formula.  This is obtained by showing that the number of vertices in a Galton-Watson tree whose out-degree lies in a given set is distributed like the number of vertices in a Galton-Watson tree with a related offspring distribution.
\end{abstract}

\maketitle

\section{Introduction}
Recently there has been considerable interest in the literature in studying the asymptotic properties of random trees.  One of the first results that obtained rescaled convergence of trees themselves (as opposed to convergence of some statistics of the trees) was a result of Aldous in \cite{Aldo93}, which proved that critical Galton-Watson trees with finite variance offspring distribution conditioned on their number of vertices converge to the Brownian continuum random tree when properly rescaled.  Subsequent work on scaling limits of random trees proceeded largely in two directions, one focusing on conditioned Galton-Watson trees and other trees with nice encodings by continuous functions see \cite{Gall06} for an overview) and another that focused on Markov branching trees.  While Galton-Watson trees conditioned on their number of vertices are Markov branching trees, the work on general Markov branching trees typically includes consistency hypotheses that are not satisfied by conditioned Galton-Watson trees  (see e.g. \cite{HMPW08}).  More recently, techniques for handling Markov branching trees that don't necessarily satisfy consistency relations and whose size is their number of leaves (or vertices) were developed in \cite{HaMi10}.  We will generalize this approach to Markov branching trees whose size is their number of vertices whose out-degree falls in a given set.  Our result for Markov branching trees, which requires a bit of setup to state rigorously, is Theorem \ref{th general A}.

As a consequence of Theorem \ref{th general A} we obtain a new theorem for scaling limits of Galton-Watson trees.  In the literature on Galton-Watson trees conditioned on their total number of vertices the techniques employed generally rely heavily on the precise form of the conditioning and cannot easily be modified.  By using the framework of Markov branching trees instead, we are able to modify the conditioning to condition on the number of vertices with out-degree in a given set.  In particular, we prove the following theorem, the notation for which will be fully defined later.

\begin{theorem}\label{th main result}
Let $T$ be a critical Galton-Watson tree with offspring distribution $\xi$ such that $0<\sigma^2=\var(\xi)<\infty$ and fix $A\subseteq \{0,1,2,\dots\}$.  Suppose that for sufficiently large $n$ the probability that $T$ has exactly $n$ vertices with out-degree in $A$ is positive, and for such $n$ let $T^A_n$ be $T$ conditioned to have exactly $n$ vertices with out-degree in $A$, considered as a rooted unordered tree with edge lengths $1$ and the uniform probability distribution $\mu_{\partial_A T^A_n}$ on its vertices with out-degree in $A$.  Then
\[\frac{1}{\sqrt{n}} T^A_n \overset{d}{\to} \frac{2}{\sigma\sqrt{\xi(A)}} T_{1/2,\nu_2},\]
where the convergence is with respect to the rooted Gromov-Hausdorff-Prokhorov topology and $T_{1/2,\nu_2}$ is the Brownian continuum random tree.
\end{theorem}

In the case $A=\Z^+= \{0,1,2,\dots\}$ we recover the classical result about the scaling limit of a Galton-Watson tree conditioned on its number of vertices first obtained in \cite{Aldo93}.  For other choices of $A$ the result appears to be new.  We note, however, that subsequent to the appearance of a draft of this article on the arXiv (\cite{Rizz11}), Kortchemski obtained similar results by different techniques in \cite{Kort11}.  The condition that for sufficiently large $n$ the probability that $T$ has exactly $n$ vertices with out-degree in $A$ is positive is purely technical and could be dispensed with at the cost of chasing periodicity considerations through our computations.  In addition to generalizing the results of \cite{HaMi10}, the key to proving this theorem is a generalization of the classical Otter-Dwass formula, which we prove in Section \ref{subsec Otter-Dwass}.  The Otter-Dwass formula (see \cite{Pitm06}) has been an essential tool in several proofs that the Brownian continuum random tree is the scaling limit of Galton-Watson trees conditioned on their number of vertices, including the original proof in \cite{Aldo93} as well as newer proofs in \cite{LeGa10} and \cite{HaMi10}.  While we follow the approach in \cite{HaMi10}, our generalization of the Otter-Dwass formula should allow for proofs along the lines of \cite{Aldo93} and \cite{LeGa10} as well.  Furthermore, with our results here, it should be straightforward to prove the analogous theorem in the infinite variance case using the approach in \cite{HaMi10}.

This paper is organized as follows.  Section \ref{sec tree models} introduces our basic notation, as well as the Markov branching trees and continuum trees that will arise for us as scaling limits.  It concludes with our generalization of the scaling limits in \cite{HaMi10}.  Section \ref{sec gw-convergence} is devoted to our study of Galton-Watson trees.  We begin by proving our generalization of the Otter-Dwass formula and we then use this to analyze the asymptotics of the partition at the root of a Galton-Watson tree.  Bringing this all together, we finish with the proof of Theorem \ref{th main result}.

\section{Models of trees}\label{sec tree models}

\subsection{Basic notation}
Fix a countably infinite set $S$; we will consider the vertex sets of all graphs discussed to be subsets of $S$.  A rooted ordered tree is a finite acyclic graph $t$ with a distinguished vertex called the root and such that, if $v$ is a vertex of $t$, the set of vertices in $t$ that are both adjacent to $v$ and further from the root than $v$ with respect to the graph metric is linearly ordered.  Let $\T^S$ be the set of all rooted ordered trees whose vertex set is a subset of $S$.  For $t,s\in \T^S$, define $t\sim s$ if and only if there is a root and order preserving isomorphism from $t$ to $s$ and let $\T=\T^S/\sim$ be the set of rooted ordered trees considered up to root and order preserving isomorphism.  By a similar construction, we let $\T^u$ be the set of rooted unordered trees considered up to root preserving isomorphism.  If $t$ is in $\T$ or $\T^u$ and $v\in t$ is a vertex, the out-degree of $v$ is the number of vertices in $t$ that are both adjacent to $v$ and further from the root than $v$ with respect to the graph metric.  The out-degree of $v$ will simply be denoted by $\deg(v)$, since we will only ever discuss out-degrees.  Fix a set $A\subseteq \Z^+$ and for $t$ in $\T$ or $\T^u$ define $\#_A t$ to be the number of vertices in $t$ whose out-degree is in $A$.  Furthermore, we define $\T_{A,n}$ and $\T^u_{A,n}$ by
\[\T_{A,n} = \{t\in \T : \#_A t= n\} \quad\textrm{and}\quad \T^u_{A,n} = \{t\in \T^u : \#_A t= n\}.\]

\subsection{Markov branching trees}\label{subsec markov branching}
In this section we extend the notion of Markov branching trees developed in \cite{HaMi10}, where Markov branching trees were constructed separately in the cases $A=\{0\}$ and $A=\Z^+$.  Here we give a construction for general $A$.  For $n\geq 0$, let $\bar\calP_n$ be the set of partitions of $n$; that is
\[\bar\calP_n = \left\{ (\lambda_1,\dots,\lambda_p) \in \bigcup_{k=0}^\infty\{0,1,2,\dots\}^k \ : \ \lambda_1\geq \lambda _2  \geq \cdots \geq \lambda_p \textrm{ and } \sum_{k=1}^p \lambda_k=n \right\}.\]
For $\lambda\in \bar\calP_n$, let $p(\lambda)$ be the length of $\lambda$ and $m_j(\lambda)$ the number of blocks in $\lambda$ equal to $j$.  If $\lambda=(\lambda_1,\dots,\lambda_p) \in \bar\calP_n$ and $k\geq 0$ we denote by $(\lambda, 0_k)$ the partition
\[ (\lambda_1,\dots,\lambda_p,\overbrace{0,\dots, 0}^{k \ 0\textrm{'s}})\]
By convention, include $\emptyset=(0_0)$ in $\bar\calP_0$ and take $p(\emptyset)=0$.  Define $\bar\calP^A_0=\bar\calP_0$ and for $n\geq 1$, define $\bar\calP^A_n$ by
\[ \bar\calP^A_n = \{\lambda \in \bar\calP_n : p(\lambda)\notin A\} \cup  \{\lambda \in \bar\calP_{n-1} : p(\lambda)\in A\}.\]

Let $(n_k)$ be an increasing sequence of integers and let $(q_{n_k})_{k\geq 1}$ be a sequence such that for each $k$ such that $n_k\geq 1$,  $q_{n_k}$ is a probability measure on $\bar\calP^A_{n_k}$ and if $n_1=0$, $q_{0}$ is a probability measure on $\T^u_{A,0}$.  For convenience, we will often consider $q_{n_k}$ to be a measure on $\cup_{j\geq 0} \bar\calP_k$ by extending it to be $0$ off of $\bar\calP_{n_k}$. 

Our goal is to construct a sequence of laws $(\bP^q_{n_k})_{k=1}^\infty$ such that $\bP^q_{n_k}$ is a law on $\T^u_{A,n_k}$ and such that the subtrees above a vertex are conditionally independent given their sizes and the degree of that vertex.  To facilitate this, we will assume that our sequences $(q_{n_k})_{k\geq 1}$ satisfy the following condition.

\begin{condition}
\begin{enumerate}
\item[(i)] Either $n_1$ or $n_2$ is equal to $1$.
\item[(ii)] For each $k$, $q_{n_k}$ is concentrated on partitions $\lambda = (\lambda_1,\dots, \lambda_p)$ such that $q_{\lambda_i}$ is defined for all $i$.
\item[(iii)] For all $n$ such that $q_n$ is defined, $\displaystyle \sum_{k=0}^\infty q_n((n,0_k)) < 1$.
\end{enumerate} 
\end{condition}

\begin{remark} Note that these assumptions put nontrivial restrictions on the compatible sequences we consider.  For example, if $1\notin A$ and $2\in A$ then $q_2$ cannot exist because $\bar\calP^A_2 = \{(2)\}$ and we are supposing that $q_2((2))<1$.  In terms of the trees we are considering this is to be expected because if a tree has root degree $2$ then it has at least $2$ leaves and, as a result, cannot possibly have exactly $2$ vertices with out-degree in $A$.
\end{remark}

If $n_1=0$, we let $\bP^q_0=q_0$.  For $n_k\geq 1$, define $\alpha_{n_k} = \sum_{j\geq 0} q_{n_k}((n_k,0_j))$ and let $G$ be geometrically distributed with
\[ \P(G=j) = (1-\alpha_{n_k})\alpha_{n_k}^j.\] 
Construct a path with $G$ edges from a from a root to a leaf.  For each non-leaf vertex $v$ on this path, independently of everything else, construct a random variable $V$ such that $\P(V=j) = q_{n_k}((n_k,0_j))/\alpha_{n_k}$ and a vector $(T_1,\dots, T_V)$ of trees that, conditionally given $V$, are i.i.d with common distribution $\bP^q_0$.  Then attach the root of each tree $T_i$ to $v$ by an edge.  For the leaf, choose a partition $\Lambda$ according to $q_{n_k}(\cdots | \bar \calP^A_{n_k}\setminus \cup_{j\geq 0}\{(n_k,0_j)\})$ and let $(T^1_1,\dots,T^1_{p(\Lambda)})$ be a vector, independent of everything except $\Lambda$, such that conditionally given $\Lambda$ the coordinates are independent and $T^1_i$ has distribution $\bP^q_{\Lambda_i}$.  Attach the roots of these trees to the leaf by an edge to obtain a tree $T$.  We define $\bP^q_{n_k}$ to be the law of $T$.

\begin{remark}
The case when $0\in A$ and $n=1$ requires some additional comment.  In this case, there is no $\bP^q_0$, and the construction above should be interpreted as producing a path of length $G$ from a root to a leaf for $\bP^q_1$.
\end{remark}

To connect with \cite{HaMi10}, if $(n_k) = (1,2,3,\dots)$, the case $A=\{0\}$ corresponds to the $\bP^q_n$ defined in \cite{HaMi10} and the case $A=\Z^+$ corresponds to the $\mathbf{Q}^q_n$ defined in \cite{HaMi10}.  Other choices of $A$ interpolate between these two extremes.  A sequence $(T_{n_k})_{k \geq 1}$ such that for each $k$, $T_{n_k}$ has law $\bP^q_{n_k}$ for some choice of $A$ and $q$ (independent of $n$) is called a Markov branching family.  For ease of notation, we will generally drop the subscript $k$ and it will be implicit that we are only considering $n$ for which the quantities discussed are defined.

\subsection{Trees as metric measure spaces}
The trees under consideration can naturally be considered as metric spaces with the graph metric.  That is, the distance between two vertices is the number of edges on the path connecting them.  Let $(t,d)$ be a tree equipped with the graph metric.  For $a>0$, we define $at$ to be the metric space $(t,ad)$, i.e. the metric is scaled by $a$.  This is equivalent to saying the edges have length $a$ rather than length $1$ in the definition of the graph metric.  More, generally we can attach a positive length to each edge in $t$ and use these in the definition of the graph metric.  Moreover, the trees we are dealing with are rooted so we consider $(t,d)$ as a pointed metric space with the root as the distinguished point.  Moreover, we are concerned with the vertices whose out-degree is in $A$, so we attach a measure $\mu_{\partial_A t}$, which is the uniform probability measure on $\partial_A t = \{v\in t : \deg(v)\in A\}$.  If we have a random tree $T$, this gives rise to a random pointed metric measure space $(T,d,\textrm{root},\mu_{\partial_AT})$.  To make this last concept rigorous, we need to put a topology on pointed metric measure spaces.  This is hard to do in general, but note that the pointed metric measure spaces that come from the trees we are discussing are compact.

Let $\M$ be the set of equivalence classes of compact pointed metric measure spaces (equivalence here being up to isometries that preserve the distinguished point and the measure).  It is worth pointing out that $\M$ actually is a set in the sense of Zermelo-Fraenkel set theory with the axiom of choice, though this takes some work to show (more precisely, there exists a set $\M$ of compact pointed metric measure spaces such that every compact pointed metric measure space is equivalent to exactly one element of $\M$).  We metrize $\M$ with the pointed Gromov-Hausdorff-Prokhorov metric (see \cite{HaMi10}).  Fix $(X,d,\rho,\mu), (X',d',\rho,\mu') \in \M$ and define
\[ d_{\textrm{GHP}}(X,X') = \inf_{(M,\delta)} \inf_{\phi:X\to M \atop \phi' :X'\to M}\left[ \delta(\phi(\rho),\phi'(\rho')) \vee \delta_H(\phi(X),\phi'(X'))\vee \delta_P(\phi_*\mu,\phi'_*\mu')\right],\]
where the first infimum is over metric spaces $(M,\delta)$, the second infimum if over isometric embeddings $\phi$ and $\phi'$ of $X$ and $X'$ into $M$, $\delta_H$ is the Hausdorff distance on compact subsets of $M$, and $\delta_P$ is the Prokhorov distance between the pushforward $\phi_*\mu$ of $\mu$ by $\phi$ and the pushforward $\phi'_*\mu'$ of $\mu'$ by $\phi'$.  Again, the definition of this metric has potential to run into set-theoretic difficulties, but they are not terribly difficult to resolve.

\begin{proposition}[Proposition 1 in \cite{HaMi10}] 
The space $(\M, d_{\textrm{GHP}})$ is a complete separable metric space.
\end{proposition}

We will not need many technical facts about Gromov-Hausdorff-Prokhorov convergence, but we will need the following lemma that shows how the pointed Gromov-Hausdorff-Prokhorov metric can be controlled in some cases by the pointed Hausdorff metric.

\begin{lemma}\label{lemma ghp control}
Let $(E,d)$ be a metric space and let $U, V \subseteq E$ be compact subsets, each with a distinguished point.  Let $\mu=n^{-1}\sum_{i=1}^n \delta_{x_i}$ and $\nu=n^{-1}\sum_{i=1}^n \delta_{y_i}$ be probability measures whose supports are contained in $U$ and $V$ respectively. Let $f:\{x_1,\dots,x_n\}  \to \{y_1,\dots, y_n\}$ be a bijection.  Then
\[d_{\textrm{GHP}}((U,\mu),(V,\nu)) \leq d_{\textrm{H}}(U,V)+\max\{d(x_i,f(x_i)), 1\leq i\leq n\}, \]
where $d_{\textrm{H}}$ is the pointed Hausdorff metric on pointed subsets of $E$.
\end{lemma}

\begin{proof}
Let $B$ be a measurable set  and define $B^\epsilon=\{y\in E \ : \ d(y,B)<\epsilon\}$.  Suppose that $\epsilon > \max\{d(x_i,f(x_i)), 1\leq i\leq n\}$.  Observe that if $x_i \in B$ then $f(x_i) \in B^\epsilon$.  Since $f$ is a bijection, we see that $|B\cap \{x_1,\dots, x_n\}| \leq |B^\epsilon \cap \{y_1,\dots, y_n\}|$, so $\mu(B)\leq \nu(B^\epsilon)$.  The reverse inequality is shown similarly and the conclusion follows.
\end{proof}

An $\R$-tree is a complete metric space $(T,d)$ with the following properties:
\begin{itemize}
\item For $v,w\in T$, there exists a unique isometry $\phi_{v,w}:[0,d(v,w)]\to T$ with $\phi_{v,w}(0)=v$ to $\phi_{v,w}(d(v,w))=w$.
\item For every continuous injective function $c:[0,1]\to T$ such that $c(0)=v$ and $c(1)=w$, we have $c([0,1]) = \phi_{v,w}([0,d(v,w)])$.
\end{itemize}

If $(T,d)$ is a compact $\R$-tree, every choice of root $\rho\in T$ and probability measure $\mu$ on $T$ yields an element $(T,d,\rho,\mu)$ of $\M$.  With this choice of root also comes a height function $\Ht(v) = d(v,\rho)$.  The leaves of $T$ can then be defined as a point $v\in T$ such that $v$ is not in $[[\rho,w[[ \deq \phi_{\rho,w}([0,\Ht(w)))$ for any $w\in T$.  The set of leaves is denoted $\mathcal{L}(T)$.  

\begin{definition}
A continuum tree is an $\R$-tree $(T,d,\rho,\mu)$ with a choice of root and probability measure such that $\mu$ is non-atomic, $\mu(\mathcal{L}(T))=1$, and for every non-leaf vertex $w$, $\mu\{v\in T : [[\rho,v]]\cap [[\rho,w]] = [[\rho,w]]\}>0$.
\end{definition}

The last condition says that there is a positive mass of leaves above every non-leaf vertex.  We will usually just refer to a continuum tree $T$, leaving the metric, root, and measure as implicit.  A continuum random tree (CRT) is an $(\M,d_{GHP})$ valued random variable that is almost surely a continuum tree.  The continuum random trees we will be interested in are those associated with self-similar fragmentation processes.    

\subsection{Self-similar fragmentations}
In this section, we recall the basic facts about self-similar fragmentations, see \cite[Chapter 3]{Bert06} for details.  For any set $B$, let $\calP_B$ be the set of countable partitions of $B$, i.e. countable collections of disjoint sets whose union is $B$.  For $n\in \overline{\N}\deq \N\cup \{\infty\}$, let $\calP_n \deq \calP_{[n]}$, where $[n]=\{1,\dots,n\}$ and $[\infty]=\N$.  Suppose that $\pi = (\pi_1,\pi_2,\dots)\in \calP_n$ (here and throughout we index the blocks of $\pi$ in increasing order of their least elements), and $B\subseteq \N$.  Define the restriction of $\pi$ to $B$, denoted by $\pi_{|B}$ or $\pi\cap B$, to be the partition of $[n]\cap B$ whose elements are the blocks $\pi_i\cap B$, $i\geq 1$.  We equip $\calP_n$ with the topology induced by the metric $d(\pi,\sigma) \deq (\inf\{i \ : \pi\cap [i]\neq \sigma\cap [i]\})^{-1}$.

\begin{definition}[Definition 3.1 in \cite{Bert06}]
Consider two blocks $B\subseteq B'\subseteq \N$.  Let $\pi$ be a partition of $B$ with $\#\pi=n$ non-empty blocks ($n=\infty$ is allowed), and $\pi^{(\cdot)}=\{\pi^{(i)}, i\in [n]\}$ be a sequence in $\calP_{B'}$.  For every integer $i$, we consider the partition of the $i$-th block $\pi_i$ of $\pi$ induced by the $i$-th term $\pi^{(i)}$ of the sequence $\pi^{(\cdot)}$, that is,
\[\pi^{(i)}_{|\pi_i}  = \left(\pi^{(i)}_j\cap\pi_i\: j\in\N\right).\]
As $i$ varies in $[n]$, the collection $\left\{\pi^{(i)}_j\cap\pi_i: i,j\in \N\right\}$ forms a partition of $B$, which we denote by $\frag(\pi,\pi^{(\cdot)})$ and call the fragmentation of $\pi$ by $\pi^{(\cdot)}$.
\end{definition}

We will use the $\frag$ function to define the transition kernels of our fragmentation processes.  Define
\[\S^\downarrow = \left\{(s_1,s_2,\dots) : s_1\geq s_2\geq\cdots\geq 0, \sum_{i=1}^\infty s_i \leq 1\right\},\]
and 
\[\S_1 = \left\{(s_1,s_2,\dots)\in [0,1]^\N \ |\ \sum_{i=1}^\infty s_i\leq 1 \right\},\]
and endow both with the topology they inherit as subsets of $[0,1]^\N$ with the product topology.  Observe that $\S^\downarrow$ and $\S_1$ are compact.  For a partition $\pi \in \calP_\infty$, we define the asymptotic frequency $|\pi_i|$ of the $i$'th block by $|\pi_i| = \lim_{n\to\infty} n^{-1}|\pi_i\cap[n]|$, provided this limit exists.  If all of the blocks of $\pi$ have asymptotic frequencies, we define $|\pi|\in \S_1$ by $|\pi| = (|\pi_1|,|\pi_2|,\dots)$.

\begin{definition}[Definition 3.3 in \cite{Bert06}]
Let $(\Pi(t),t\geq 0)$ be an exchangeable, c\`adl\`ag $\calP_\infty$-valued process such that $\Pi(0)=\mathbf{1}_\N\deq(\N,0,\dots)$ and such that 
\begin{enumerate}
\item $\Pi(t)$ almost surely possesses asymptotic frequencies $|\Pi(t)|$ simultaneously for all $t\geq 0$ and
\item if we denote by $B_i(t)$ the block of $\Pi(t)$ which contains $i$, then the process $t\mapsto |B_i(t)|$ has right-continuous paths.
\end{enumerate}
We call $\Pi$ a self-similar fragmentation process with index $\alpha\in \R$ if and only if, for every $t,t'\geq 0$, the conditional distribution of $\Pi(t+t')$ given $\F_t$ is that of the law of $\frag(\pi,\pi^{(\cdot)})$, where $\pi=\Pi(t)$ and $\pi^{(\cdot)}=\left(\pi^{(i)},i\in \N\right)$ is a family of independent random partitions such that for $i\in \N$, $\pi^{(i)}$ has the same distribution as $\Pi(t'|\pi_i|^\alpha)$.
\end{definition}

Suppose, for the moment, that $\Pi$ is a self-similar fragmentation process with $\alpha=0$ (these are also called homogeneous fragmentations).  In this case $\Pi$ is a Feller process as is $\Pi_{|[n]}$ for every $n$.  For $\pi\in \calP_n\setminus\{\mathbf{1}_{[n]}\}$, let
\[q_\pi = \Lim{t\to 0^+} \frac{1}{t} \P(\Pi_{|[n]}(t)=\pi).\]
For $\pi\in \calP_n$ and $n'\in \{n,n+1,\dots,\infty\}$, define $\calP_{n',\pi} = \{\pi'\in\calP_{n'} : \pi'_{|[n]}=\pi\}$.
One way to construct a homogeneous fragmentation is by specifying these transition rates, which can be done using measures on $\S^\downarrow$.  The basic building blocks are the so called paintbox partitions that correspond to elements of $\S^\downarrow$ as follows.  Let $(U_1,U_2,\dots)$ be a sequence of i.i.d uniform $(0,1)$ random variables and for $\mathbf{s}\in \S^\downarrow$, let $\rho_{\mathbf{s}}$ be the law of the partition defined by $i$ and $j$ are in the same block if and only if there exists $k$ such that for $\ell \in \{i,j\}$
\[\sum_{m=1}^ks_m \leq U_\ell <\sum_{m=1}^{k+1}s_m.\]

\begin{theorem}[see Propositions 3.2 and 3.3 in \cite{Bert06}] 
Given a measure $\nu$ on $\S^\downarrow$ such that $\nu(\{\mathbf{1}\})=0$ and $\int_{\S^{\downarrow}}(1-s_1)\nu(d\mathbf{s})<\infty$, define a measure $\rho_\nu$ on $\calP_\infty$ by
\[\rho_\nu(\cdot) = \int_{\mathbf{s}\in\S^{\downarrow}} \rho_{\mathbf{s}}(\cdot) \nu(d\mathbf{s}).\]
Then there is a unique homogeneous fragmentation such that $q_\pi = \rho_\nu(\calP_{\infty,\pi})$ for all $\pi\in \cup_n (\calP_n\setminus\{\mathbf{1}_{[n]}\})$.
\end{theorem}
If $\Pi^0(t)$ is a homogenous fragmentation corresponding to a measure $\nu$ as in the above theorem, we will call $\nu$ the splitting measure $\Pi^0(t)$.  If $\alpha<0$, we can construct an $\alpha$-self-similar fragmentation from $\Pi^0$ by a time change.  
Let $\pi(i,t)$ be the block of $\Pi^0$ that contains $i$ at time $t$ and define
\[ T_i(t) = \inf\left\{u\geq 0 : \int_0^u|\pi(i,r)|^{-\alpha} dr>t \right\}.\]
For $t\geq 0$, let $\Pi(t)$ be the partition such that $i,j$ are in the same block of $\Pi(t)$ if and only if they are in the same block of $\Pi^0(T_i(t))$.  Then $(\Pi(t),t\geq 0)$ is a self-similar fragmentation.   Moreover, for we call $(\alpha,\nu)$ a fragmentation pair and and we call  $(\Pi(t),t\geq 0)$ a fragmentation with characteristics $(\alpha,\nu)$.

We will need trees associated to fragmentations with characteristics $(\alpha,\nu)$, where $\alpha<0$ and $\nu(\sum_is_i<1)=0$; henceforth, we let $\Pi$ be such a self-similar fragmentation.  The tree associated a fragmentation processes $\Pi$ is a continuum random tree that keeps track of when blocks split apart and the sizes of the resulting blocks.  For a continuum tree $(T,\mu)$ and $t\geq 0$, let $T_1(t),T_2(t),\dots$ be the tree components of $\{v\in T : \mathrm{ht}(v)>t\}$, ranked in decreasing order of $\mu$-mass.  We call $(T,\mu)$ self-similar with index $\alpha<0$ if for every $t\geq 0$, conditionally on $(\mu(T_i(t)),i\geq 1)$, $(T_i(t),i\geq 1)$, equipped with the restriction of $\mu$ normalized to be a probability measure, has the same law as $(\mu(T_i(t))^{-\alpha}T^{(i)},i\geq 1)$ where the $T^{(i)}$'s are independent copies of $T$. 

The following summarizes the parts of Theorem 1 and Lemma 5 in \cite{HaMi04} that we will need.

\begin{theorem}\label{theorem existence of trees}
Let $\Pi$ be a $(\alpha,\nu)$-self-similar fragmentation with $\alpha<0$ and $\nu$ as above and let $F\deq |\Pi|^\downarrow$ be its ranked sequence of asymptotic frequencies.  There exists an $\alpha$-self-similar CRT $(T_{-\alpha,\nu},\mu_{-\alpha,\nu})$ such that, writing $F'(t)$ for the decreasing sequence of masses of the connected components of $\{v\in T_{-\alpha,\nu}:ht(v)>t\}$, the process $(F'(t),t\geq 0)$ has the same law as $F$.  Furthermore, $T_{-\alpha,\nu}$ is a.s. compact.
\end{theorem}  

The choice of where to put negative signs in the notation in the above theorem is to conform with the notation of \cite{HaMi10}.  

\begin{definition} The Brownian CRT is the $-1/2$-self-similar random tree with dislocation measure $\nu_2$ given by
\[\int_{\S^\downarrow} \nu_2(d\mathbf{s})f(\mathbf{s})= \int_{1/2}^1\sqrt{\frac{2}{\pi s_1^3(1-s_1)^3}} ds_1 f(s_1,1-s_1,0,0,\dots).\]
\end{definition}

\subsection{Convergence of Markov branching trees}
We first recall some of the main results of \cite{HaMi10}.  Let $A\subseteq \Z^+$ and let $(q_n)$ be a compatible sequence of probability measures satisfying the conditions of Section \ref{subsec markov branching}.  Define $\bar q_n$ to be the push forward of $q_n$ to $\S^\downarrow$ by $\lambda \mapsto \lambda/\sum_i\lambda_i$.

\begin{theorem}[Theorems 1 and 2 in \cite{HaMi10}] \label{th main HaMi10}
Suppose that $A=\{0\}$ or $A=\Z^+$.  Further suppose that there is a fragmentation pair $(-\gamma,\nu)$ with $0<\gamma < 1$ and a function $\ell:(0,\infty)\to (0,\infty)$, slowly varying at $\infty$ (or $\gamma=1$ and $\ell(n)\to 0$) such that, in the sense of weak convergence of finite measures on $\S^\downarrow$, we have
\[ n^\gamma \ell(n)(1-s_1)\bar q_n(d\mathbf{s}) \rightarrow (1-s_1)\nu(d\mathbf{s}).\]
Let $T_n$ have law $\bP^q_n$ and view $T_n$ as a random element of $\M$ with the graph distance and the uniform probability measure $\mu_{\partial_A T_n}$ on $\partial_A T_n = \{v\in T_n : \deg\ v\in A\}$.
Then we have the convergence in distribution
\[ \frac{1}{n^\gamma \ell(n)} T_n \rightarrow T_{\gamma,\nu},\]
with respect to the rooted Gromov-Hausdorff-Prokhorov topology.
\end{theorem}

The case where $A=\{0\}$ this is a special case of Theorem 1 in \cite{HaMi10} and the case $A=\Z^+$ is Theorem 2 in the same paper.  The case $A=\Z^+$ is proved by reduction to the $A=\{0\}$ case.  We extend this to the case of general $A$ containing $0$.  Fix, for the moment, a tree $t$.  For a vertex $v\in t$, let $t_v$ be the subtree of $t$ above $v$.  We call the elements of $\{t_v : v\in t\}$ the fringe subtrees of $t$ and say that $v$ is the root of $t_v$.  Furthermore, we partially order the fringe subtrees of $t$ by inclusion.

\begin{theorem}\label{th general A}
The conclusions of Theorem \ref{th main HaMi10} are valid if the only assumption on $A\subseteq \Z^+$ is that $0\in A$.  If $0 \notin A$, the conclusions of Theorem \ref{th main HaMi10} remain valid if we further assume that
\[d_{\textrm{GHP}}\left(\frac{1}{n^\gamma \ell(n)} T_n, \frac{1}{n^\gamma \ell(n)} \tilde T_n\right) \to 0\]
in probability, where $ \tilde T_n$ is the tree obtained from $T_n$ by deleting all maximal fringe subtrees of $T_n$ that contain no vertices with out-degree in $A$.
\end{theorem}

\begin{remark}
The condition that
\[d_{\textrm{GHP}}\left(\frac{1}{n^\gamma \ell(n)} T_n, \frac{1}{n^\gamma \ell(n)} \tilde T_n\right) \to 0\]
in probability expresses a relationship between the number of subtrees that have no vertices with out-degree in $A$ and the height of these trees.  For example, it is always satisfied if there is a uniform upper bound on the trees in the support of $\bP^q_0$.  In our application to Galton-Watson trees, this condition is easily checked, but it may be of interest to investigate general conditions on $\bP^q_0$ and $(q_{n_k})_{k\geq 1}$ that imply this convergence.
\end{remark} 

We will prove Theorem \ref{th general A} by reduction to Theorem \label{th main HaMi10}.  For a tree $t$,  we define $t^\circ$ to be the tree obtained from $t$ as follows:
\begin{enumerate}
\item Attach a mark to each vertex in $t$ whose out-degree is in $A$.
\item For each vertex $v$, delete the tree $t_v$ (which is the subtree of $t$ above $v$) if it contains no marked vertices.
\item For each vertex $v$, if $t_v$ is a maximal among fringe subtrees containing exactly one marked vertex, replace $t_v$ by a marked leaf.
\item Attach a leaf to each marked vertex that is not a leaf at the completion of Step (3).
\item Remove all the marks from the tree.
\end{enumerate}


Suppose that $(q_{n_k})_{k\geq 1}$ satisfies the hypotheses above.  For notational convenience we extend $q_n$ to $\cup_k \bar\calP^A_k$ by setting $q_n(\lambda)=0$ for $\lambda\notin \bar\calP^A_n$.  Define measures
\[q_n^1(\lambda) =   \sum_{k=0}^\infty q_n((\lambda, 0_k))\cf(p(\lambda)+k\notin A) \quad \textrm{and}\quad  q^2_n(\lambda) =\sum_{k=0}^\infty q_n((\lambda,0_k))\cf(p(\lambda)+k\in A). \]
We can then define probability measures $q^\circ_n$ on $\bar\calP_n$ by $q^\circ_1(\emptyset)=1$ and for $n\geq 2$
\[ q^\circ_n(\lambda) = \begin{cases}q^1_n(\lambda) & \textrm{if } \lambda_{p(\lambda)}\geq 2 \\ q_n^1(\lambda) + q_n^2(\lambda')& \textrm{if } \lambda=(\lambda',1)\\
0 &\textrm{otherwise}.\end{cases}\]


\begin{lemma}
If $n\geq 1$ and $T_n$ is distributed like $\bP^q_n$ then $T^\circ_n$ is distributed like $\bP^{q^\circ}_n$.
\end{lemma}

\begin{proof}
We proceed by induction on $n$.    For $n=1$, it is clear from Step (3) in the construction of $T^\circ_1$ that $T^\circ_1$ is a tree with a single vertex; $\bP^{q^\circ}_n$ is concentrated on this tree since $q^\circ_1(\emptyset)=1$.  For $n\geq 2$, we condition on the root partition.  First suppose that the root partition is not of the form $(n,0_k)$.  In this case, in both $T_n^\circ$ and a tree with law $\bP^{q^\circ}_n$, given the root partition the subtrees attached to the root are independent  and have the appropriate distributions by the inductive hypothesis.  When the root partition is of the form $(n,0_k)$ note that in both trees the lowest vertex whose out-degree is not $1$ is attached to the root by a path of geometric length and the parameters of the geometric variable are the same.  Moreover, in both $T_n^\circ$ and a tree distributed like $\bP^{q^\circ}_n$, the fringe subtree above lowest vertex whose out-degree is not $1$ is distributed like the tree conditioned on the root vertex not being of the form $(n,0_k)$.  Thus, by induction, we need only check that the laws of the partitions at the root agree.  This, however, is immediate from the construction of $q_n^\circ$.  
\end{proof}

\begin{lemma}\label{le follower}
If 
\[n^\gamma \ell(n)(1-s_1)\bar q_n(d\mathbf{s}) \rightarrow (1-s_1)\nu(d\mathbf{s}),\]
then
\[ n^\gamma \ell(n)(1-s_1)\bar q^\circ_n(d\mathbf{s}) \rightarrow (1-s_1)\nu(d\mathbf{s}).\]
\end{lemma}

\begin{proof}
Let $f:\S^\downarrow\to \R$ be Lipschitz continuous (with respect to the uniform norm) with both the uniform norm and Lipschitz constant bounded by $K$.  For $\lambda=(\lambda_1,\dots,\lambda_r,0_k)$ with $\lambda_r\geq 1$, define $r(\lambda)=r$ and $\iota(\lambda) = (\lambda_1,\dots,\lambda_r,1,0_k)$.  Observe that for $\lambda\in \bar\calP_n$,
\[\left| f\left(\frac{\iota(\lambda)}{n+1}\right) - f\left(\frac{\lambda}{n}\right)\right| \leq K\sum_{i=1}^{r(\lambda)} \frac{\lambda_i}{n(n+1)} + \frac{K}{n+1} = \frac{2K}{n+1}.\]
Letting $g(\mathbf{s}) = (1-s_1)f(\mathbf{s})$, we have
\[\begin{split} |\bar q_n^\circ(g)-\bar q_n(g)|& \leq \sum_{\lambda\in \bar\calP^A_n\cap \bar\calP_{n-1}} q^2_n(\lambda)\left| \left(1-\frac{\lambda_1}{n}\right) f\left(\frac{\iota(\lambda)}{n}\right) - \left(1-\frac{\lambda_1}{n-1}\right)f\left(\frac{\lambda}{n-1}\right)\right|\\
&\leq \sum_{\lambda\in \bar\calP^A_n\cap \bar\calP_{n-1}} q^2_n(\lambda)\left(\frac{K\lambda_1}{n(n-1)}+\frac{2K}{n}\right) \\
& \leq \frac{3K}{n}.
\end{split}\]
Multiplying by $n^\gamma \ell(n)$, we see that this upper bound goes to $0$ and the result follows.
\end{proof}

\begin{proof}[Proof of Theorem \ref{th general A}]
Let $ T^1_n$ be the tree obtained from $T_n$ by removing all of the maximal fringe subtrees containing no vertices with out-degree in $A$.  By hypothesis, we have
\[d_{\textrm{GHP}}\left( \frac{1}{n^\gamma \ell(n)} T^1_n, \frac{1}{n^\gamma \ell(n)} T_n\right) \to 0\]
in probability.  Note that leaves in $T^1_n$ correspond to the unique vertices with out-degree in $A$ in the fringe subtrees of $T_n$ that are maximal among fringe subtrees with exactly one vertex with out-degree in $A$.  By construction, these vertices were attached to the rest of $T_n$ by independent spines of geometric length.  Let $T^2_n$ be the tree obtained from $T^1_n$ by replacing these spines by a single vertex and moving the mass from the tip of this spine to the single vertex.  Since the maximum length of the spines altered in the creation of $T^2_n$ is $O(\log(n))$ in probability, we have 
\[d_{\textrm{GHP}}\left( \frac{1}{n^\gamma \ell(n)} T^1_n, \frac{1}{n^\gamma \ell(n)} T^2_n\right) \to 0\]
in probability.  Here we have applied Lemma \ref{lemma ghp control} to obtain the Prokhorov part of this convergence with $f$ being the map that takes a vertex at the top of a spine to the single vertex replacing the spine.  Finally, observe that $T^\circ_n$ differs from $T^2_n$ only by the attachment of some additional leaves with a corresponding adjustment of the mass measure.  Thus
\[d_{\textrm{GHP}}\left( \frac{1}{n^\gamma \ell(n)} T^2_n, \frac{1}{n^\gamma \ell(n)} T^\circ_n\right) \to 0\]
in probability.  Since $(n^\gamma \ell(n))^{-1}T_n^\circ \to T_{\gamma,\nu}$ by Lemma \ref{le follower} and Theorem \ref{th main HaMi10}, $(n^\gamma \ell(n))^{-1}T_n \to T_{\gamma,\nu}$ as well.
\end{proof}

\section{Galton-Watson trees} \label{sec gw-convergence}

Let $\xi=(\xi_i)_{i\geq 0}$ be a probability distribution with mean less than or equal to $1$, and assume that $\xi_1<1$.  A Galton-Watson tree with offspring distribution $\xi$ is a random element $T$ of $\T$ with law
\[ \GW_\xi(t) = \P(T=t) = \prod_{v\in t} \xi_{\deg(v)}.\]
The fact that $\xi$ has mean less than or equal to $1$ implies that the right hand side defines an honest probability distribution of $\T$.

\subsection{Otter-Dwass type formulae} \label{subsec Otter-Dwass}
In this section we develop a transformation of rooted ordered trees that takes Galton-Watson trees to Galton-Watson trees.  This transformation is motivated by the observation that the number of leaves in a Galton-Watson tree is distributed like the progeny of a Galton-Watson tree with a related offspring distribution.  This simple observation was first made in \cite{Mina05}.  Let $\xi=(\xi_i)_{i\geq 0}$ be a probability distribution with mean less than or equal to $1$, and assume that $\xi_1<1$.  Let $T$ be a Galton-Watson tree with offspring distribution $\xi$ and let
\[C(z) = \sum_{i=1}^\infty \P(\#_{\{0\}}T = i)z^i  \]
be the probability generating function of the number of leaves of $T$.  Furthermore, let
\[ \phi(z) =\sum_{i=0}^\infty \xi_{i+1}z^i.\]
Decomposing by the root degree, we see that $C(z)$ satisfies the functional equation 
\[C(z) = \xi_0z+C(z)\phi(C(z)).\] 
Solving for $C(z)$ yields
\begin{equation} \label{fundamental} C(z)= z\left(\frac{\xi_0}{1-\phi(C(z))}\right).\end{equation}
Define
\begin{equation}\label{eq derived dist} \theta(z) = \frac{\xi_0}{1-\phi(z)}.\end{equation}
Observe that $\theta$ has nonnegative coefficients, $[z^0]\theta(z) = \xi_0/(1-\xi_1)$ and $\theta(1)=1$.  Thus the coefficients of $\theta$ are a probability distribution, call it $\zeta=(\zeta_i)_{i\geq 0}$.  

\begin{proposition}\label{prop leaves to vertices}
Let $T$ be a Galton-Watson tree with offspring distribution $\xi$ and let $T'$ be a Galton-Watson tree with offspring distribution $\zeta$ where $\xi$ and $\zeta$ are related as above.  Then for all $k\geq 1$, $\P(\#_{\{0\}}T=k)=\P(\#_{\Z^+}T'=k)$.  Also, $T'$ is critical (subcritical) if and only if $T$ is critical (subcritical).
\end{proposition}

\begin{proof}
Let $\tilde C(z)$ be the probability generating function of $\#_{\Z^+} T'$.  A similar computation as the one above for $C(z)$ shows that $\tilde C(z)$ satisfies the functional equation $\tilde C(z) = z \theta(\tilde C(z))$, which is the same functional equation we showed $C(z)$ solves in Equation \ref{fundamental}.  The Lagrange inversion formula thus implies that $C(z)$ and $\tilde C(z)$ have the same coefficients, i.e. that $\P(\#_{\{0\}}T=k)=\P(\#_{\Z^+}T'=k)$ for all $k$.  From Equation \eqref{eq derived dist},  we see that $\theta'(z) = \xi_0\phi'(z) (1-\phi(z))^2$.  Let $\alpha$ be the mean of $\xi$.  A short computation shows that that $\phi'(1) = \alpha -1 + \xi_0$.  Since $\phi(1) = 1-\xi_0$, we see that the mean of $\zeta$ is
\[ \theta'(1) = \frac{\xi_0 (\alpha -1 + \xi_0)}{\xi_0^2} = \frac{\alpha -1 + \xi_0}{\xi_0},\]
and the criticality (subcriticality) claim follows.  The same approach can also be used to obtain higher moments of $\zeta$.
\end{proof}

\begin{corollary}\label{cor Otter-Dwass}
Let $\F_n$ be an ordered forest of $n$ independent Galton-Watson trees all with offspring distribution $\xi$.  Let $\zeta$ be related to $\xi$ as in Proposition \ref{prop leaves to vertices}.  Let $(X_i)_{i\geq 1}$ be an i.i.d.\!\! sequence of $\zeta$ distributed random variables and let $S_k =\sum_{i=1}^k (X_i-1)$.  Let $\#_{\{0\}}\F_n$ denote the number of leaves in $\F_n$.  Then for $1\leq k\leq n$
\[ \P(\#_{\{0\}}\F_k=n) = \frac{k}{n} \P(S_n=-k).\]
\end{corollary} 

\begin{proof}
This follows immediately from Proposition \ref{prop leaves to vertices} and the Otter-Dwass formula (see \cite{Pitm06}).
\end{proof}

This relationship between $T$ and $T'$ can also be proved in a more probabilistic fashion.  Indeed, by taking a more probabilistic approach we can get a more general result that includes the results in \cite{Mina05} as a special case.  The idea is a lifeline construction.  You follow the depth-first walk around the tree and when you encounter a vertex whose degree is in $A$ you label any unlabeled edges on the path from this vertex to the root with the label of this vertex.  This labeled path can be considered the lifeline of the vertex.  A new tree is constructed by letting the root be first vertex encountered whose degree is in $A$ and attaching vertices whose degree is in $A$ to the earliest vertex whose lifeline touches its own.  We will now go through this construction formally.  

\begin{figure}
\[ 
\xymatrix{
  &\ \ \bullet^2 \ar@{-}[d]^2 & \bullet^3 \ar@{-}[dl]^3 \\
\ \ \bullet^1\ar@{-}[d]^1 & \bullet \ar@{-}[d]^2 &  \bullet^4 \ar@{-}[dl]^4\\
\bullet \ar@{-}[dr]^1 &  \bullet \ar@{-}[d]^2 & \ \  \bullet^5 \ar@{-}[d]^5  & & & & & \ \ \bullet^3 \ar@{-}[d] & \bullet^4 \ar@{-}[dl] \\
&\bullet \ar@{-}[dr]^1 & \bullet \ar@{-}[d]^5 & \bullet^6  \ar@{-}[dl]^6 & \ar@{|->}[rr]^\vee & & & \ \ \bullet^2 \ar@{-}[dr] & \ \  \bullet^5 \ar@{-}[d] &  \bullet^6 \ar@{-}[dl]  \\
&& \bullet & & & & & &\ \ \bullet^1
} 
\]
\caption{A colored tree and its image under $\vee$ when $A=\{0\}$}
\label{fig check map}
\end{figure}
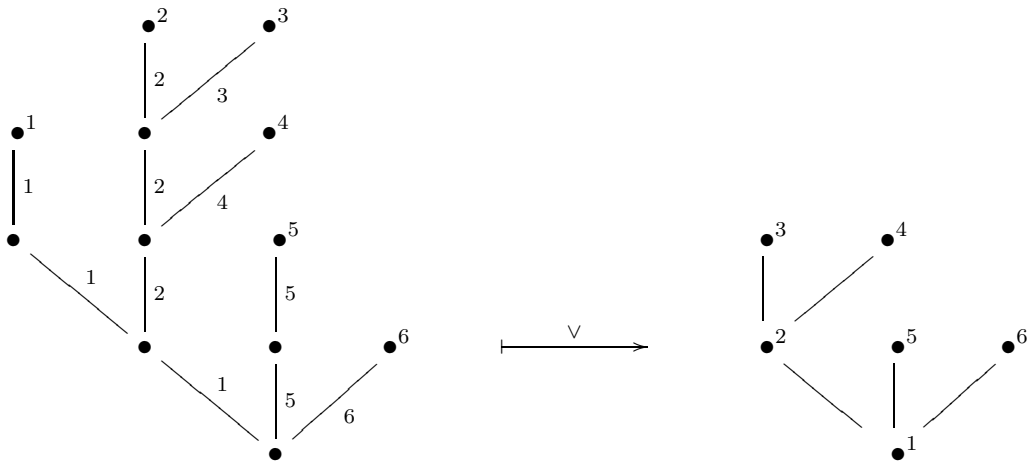

Suppose that $t\in \T_{A,n}$ for some $n\geq 1$ and label the vertices of $t$ whose out-degrees are in $A$ by the order they appear in the depth first walk of $t$.  We will now color a subset of the edges of $t$.  For each edge $e\in t$, let $t_e$ be the component of $t\setminus e$ that does not contain the root.  Color $e$ with the smallest number that is the label of a vertex in $t_e$, leaving $e$ uncolored if $t_e$ contains no labeled vertices.  Note that for any $1\leq k\leq n$ the subtree spanned by the vertices with labels $\{1,\dots,k\}$ is colored by $\{1,\dots,k\}$ and an edge is colored by an element in $\{1,\dots,k\}$ if and only if it is in this subtree.  Furthermore, the path from the vertex with label $k$ to any edge colored $k$ contains only edges colored $k$.  See Figure \ref{fig check map} for an example of such a coloring when $A=\{0\}$.  Call two edges of $t$ coincident if they share a common vertex.

\begin{lemma}\label{le basic label prop}
If $t$ is colored as above and $2\leq j\leq n$, then there is exactly one edge colored $j$ that is coincident to an edge with a smaller color.
\end{lemma}

\begin{proof}
First we show existence.  Consider the path from vertex $j$ to the root.  Let $e$ be the last edge this path that is not contained in the subtree spanned by vertices $\{1,\dots, j-1\}$.  By construction this edge is colored $j$ and is coincident to an edge colored by an element of $\{1,\dots, j-1\}$.

To see uniqueness, suppose that $f$ is an edge with the desired properties.  Then $f$ is on the path from $j$ to the root and $f$ is coincident to an edge in the subtree spanned by vertices $\{1,\dots, j-1\}$.  If $f$ contains the root, then $f$ is the last edge on the path from $j$ to the root that is colored $j$, i.e. $f=e$.  Otherwise, after $f$, we finish the path from $j$ to the root within this subtree.  Hence $f$ is the last edge on the path from $j$ to the root that is colored $j$ and again $f=e$.  
\end{proof}

With $t$ labeled as above we define a rooted plane tree with $n$ vertices, called the life-line tree and denoted $\check t$, as follows.  The vertex set of $\check t$ is $\{1,2,\dots,n\}$, $1$ is the root.  Furthermore, if $i<j$, $i$ is adjacent to $j$ if $i$ is the smallest number such that there exist coincident edges $e_1,e_2$ in $t$ with $e_1$ colored $i$ and $e_2$ colored $j$.  Finally, the children of a vertex are ordered by the appearance of the corresponding vertices with out-degree in $A$ in the depth-first search of $t$.  See Figure \ref{fig check map} for an example of this map when $A=\{0\}$.

\begin{lemma}
The life-line tree is a connected acyclic graph.
\end{lemma}

\begin{proof}
Suppose that $\check t$ has at least two components.  Let $j$ be the smallest vertex not in the same component as $1$.  By Lemma \ref{le basic label prop}, there exists $1\leq i<j$ and coincident edges $e_1,e_2$ in $t$ labeled $i$ and $j$ respectively.  Thus $i$ is adjacent to $j$, a contradiction.

Suppose that $\check t$ contains a cycle.  Let $j$ be the largest vertex in this cycle.  Then $j$ is adjacent to two smaller vertices, contradicting our definition of $\check t$.
\end{proof}

Using this construction we can prove the following theorem.

\begin{theorem}\label{theorem check transform}
Fix $A\subseteq \Z_+$ and let $T'$ be a critical or subcritical Galton-Watson tree with non-degenerate offspring distribution $\xi$ and let $T$ be $T'$ conditioned to have at least one vertex with out-degree in $A$.  Let $(X_i)_{i\geq1}$ be i.i.d random variables with distribution $\xi$.  Let $N(X) = \inf\{k : X_k\in A\}$ and $\tau_{-1}(X) = \inf\left\{k: \sum_{i=1}^k (X_i-1) = -1\right\}$.  Let $\tilde X$ be distribution like
\[\left( 1+\sum_{i=1}^{N(X)}(X_i-1)\right) \ \Big| (N(X)\leq \tau_{-1}(X)).\]
Finally, introduce the random variable $Y$ such that, conditionally given $\tilde X=k$, $Y$ is binomially distributed with parameters $k$ and $p=\P(N(X)\leq \tau_{-1}(X))$.  Then $\check T$ is distributed like a Galton-Watson tree whose offspring distribution is the law of $Y$.  Moreover, if $T'$ is critical so is $\check T$ (i.e. $\E Y = 1$) and if additionally $0<\var(\xi)=\sigma^2<\infty$ then $\var(Y)=\P(N\leq \tau_{-1})^2 \sigma^2/\xi(A)$.
\end{theorem}

The claims about the mean and variance of $Y$ in this theorem are very important for our purposes as they are what allow us to derive the constants in Theorem \ref{th main result}.  As such, we separate out the computation of these constants into its own lemma.

\begin{lemma}\label{lemma Y moments}
Let be $Y$ as defined in Theorem \ref{theorem check transform}.  If $\xi$ has mean $1$ then $\E Y=1$.  If additionally $0<\var(\xi)=\sigma^2<\infty$ then $\var(Y)=\P(N\leq \tau_{-1})^2 \sigma^2/\xi(A)$ 
\end{lemma}

\begin{proof}
Let $p = \P(N\leq\tau_{-1})$.  By conditioning on the value of $\tilde X$, we see that 
\[ \E Y =\E Y = \sum_{k=1}^\infty kp \P(\tilde X=k) = p \E \tilde X \]
and 
\[ \E Y^2 =  \sum_{k=1}^\infty (kp(1-p) +k^2p^2) \P(\tilde X=k) = p(1-p)\E \tilde X + p^2\E \tilde X ^2.\]
This reduces the problem to computing the moments of $\tilde X$.  Assume that $\xi$ has mean $1$, so that $\E X_i=1$.  By Wald's first equation we have
\[\begin{split} 1 & = \E\left[ 1+ \sum_{i=1}^{N}(X_i-1)\right]\\
& =  1+ \E\left[\sum_{i=1}^{N}(X_i-1), \ N\leq \tau_{-1}\right]  + \E\left[ \sum_{i=1}^{N}(X_i-1), \ N > \tau_{-1}\right].
\end{split} \]
Using the strong Markov property of $(X_i)_{i\geq1}$ at the stopping time $\tau_{-1}$ and Wald's first equation again we see that 
\[ \begin{split} \E\left[ \sum_{i=1}^{N}(X_i-1), \ N > \tau_{-1}\right] & = \P(N>\tau_{-1}) \left( -1+  \E\left[ \sum_{i=1}^{N}(X_i-1)\right] \right) \\
&=  -\P(N>\tau_{-1}) 
\end{split}.\]
Therefore $\E \tilde X = \P(N\leq\tau_{-1})^{-1} = p^{-1}$ and consequently $\E Y = p \E\tilde X = 1$. 

Now assume additionally that $0<\var(\xi) <\infty$.  Wald's second equation shows that 
\[\begin{split} \frac{\sigma^2}{\xi(A)}  & = \E\left(\left[ \sum_{k=1}^N (X_k-1)\right]^2\right)\\
& = \E\left(\left[ \sum_{k=1}^N (X_k-1)\right]^2, \ N\leq \tau_{-1}\right)+ \E\left(\left[ \sum_{k=1}^N (X_k-1)\right]^2, \ N>\tau_{-1}\right)
.\end{split}\]
Using the same Markov property trick as in the previous computation, we see that
\[\E\left(\left[ \sum_{k=1}^N (X_k-1)\right]^2, \ N>\tau_{-1}\right)= \P(N>\tau_{-1})\left(1 + \frac{\sigma^2}{\xi(A)}\right).\]
Consequently 
\[\begin{split} \E \tilde X^2 & =\frac{1}{\P(N\leq \tau_{-1})}  \E\left(\left[1+ \sum_{k=1}^N (X_k-1)\right]^2, \ N\leq \tau_{-1}\right)\\
& = 1 +  2p^{-1}\E\left[\sum_{i=1}^{N}(X_i-1), \ N\leq \tau_{-1}\right]  + p^{-1} \E\left(\left[ \sum_{k=1}^N (X_k-1)\right]^2, \ N\leq \tau_{-1}\right)\\
& = 1+ 2p^{-1}(1-p) + p^{-1}\left[ \frac{\sigma^2}{\xi(A)} - (1-p)\left(1 + \frac{\sigma^2}{\xi(A)}\right)\right] \\
& = p^{-1} +   \frac{\sigma^2}{\xi(A)}
.\end{split}\]
Thus
\[\var(Y)= \E Y^2-1 = (1-p) + p^2\left(p^{-1} +   \frac{\sigma^2}{\xi(A)} \right) -1 =  p^2 \frac{\sigma^2}{\xi(A)},\]
as desired.
\end{proof}

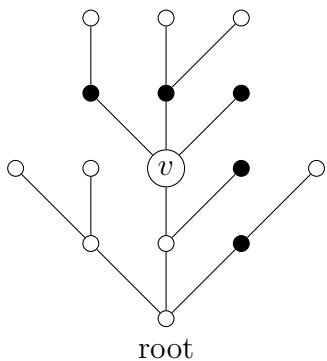
\begin{figure}[here]
\[\begin{tikzpicture}
\vertex (r) at (1,0) [label=below:root]{};
\vertex (1) at (0,1) {};
\vertex (2) at (1,1) {};
\vertex[fill] (3) at (2,1) {};
\vertex (4) at (-1,2) {};
\vertex (5) at (0,2) {};
\Vertex (6) at (1,2)  {$v$};
\vertex[fill] (7) at (0,3) {};
\vertex[fill] (8) at (1,3) {};
\vertex[fill] (9) at (2,2) {};
\vertex (10) at (3,2) {};
\vertex[fill] (11) at (2,3) {};
\vertex (12) at (0,4) {};
\vertex (14) at (2,4) {};
\vertex (13) at (1,4) {};
\path (r) edge (1) (r) edge (2) (r) edge (3) (1) edge (4) (1) edge (5) (2) edge (6) (2) edge (9) (3) edge (10) (6) edge (7) (6) edge (8) (6) edge (11) (8) edge (14) (8) edge (13) (7) edge (12);
\end{tikzpicture}\]
\caption{The black nodes are the right spine of the vertex $v$.}
\label{fig rspine}
\end{figure}

Our proof of Theorem \ref{theorem check transform} will proceed by induction.  The first step is to analyze the degree of the root of $\check T$.  To do this, we need some notation.  For a vertex $v\in t$ define the spine of $v$ to be the vertices in $t$ that are children of $v$ or children one of $v$'s ancestors.  The right spine of $v$, denoted $\textrm{rspine}(v)$, is the subset of the spine of $v$ consisting of vertices that appear after $v$ in the depth-first walk of $t$ (see Figure \ref{fig rspine}).

\begin{lemma}\label{lemma check root}
The degree of the root of $\check T$ is distributed like $Y$.
\end{lemma}

\begin{proof}
Let $(v_1,\dots, v_n)$ be the vertices of $T$ in depth first order and let $v$ be the first vertex along the depth first path of $T$ such that $\deg(v)\in A$.  The degree of the root of $\check T$ is the number of vertices in $\textrm{rspine}(v)$ that are the root of a fringe subtree with at least one vertex whose out-degree is in $A$.  Furthermore, it is clear from the definition of Galton-Watson trees that the fringe subtrees of $T$ rooted at the vertices in $\textrm{rspine}(v)$ are conditionally independent given the number of vertices in $\textrm{rspine}(v)$.  Hence, conditionally given $\#\textrm{rspine}(v)$, the degree of the root of $\check T$ has distribution $\textrm{Bin}(\#\rspine(v), \P(N\leq \tau_{-1}))$.  An easy combinatorial argument shows that if $k$ is the index of $v$ in the depth first order of $T$, then
\[\#\rspine(v) =1+ \sum_{i=1}^k (\deg(v_i)-1).\]
It is now immediate that the degree of the root of $\check T$ is distributed like $Y$.
\end{proof}

\begin{proof}[Proof of Theorem \ref{theorem check transform}]
It follows from Lemma \ref{lemma check root} that the probability that $\check T$ has a single vertex equals $\P(Y=0)$.  We will proceed by induction on the number of vertices of $t$ to show that for all $t$
\[\P(\check T = t) = \prod_{v\in t} \P(Y= \deg(v)).\]
Suppose the claim is true for all trees $s$ with fewer than $n$ vertices and let $t$ be a tree with $n$ vertices.  Let $r$ be the degree of the root of $t$ and let $(t_1,\dots,t_r)$ be the subtrees of $t$ attached to the root.

Let us consider the set $S = \{ s\in \T : \check s = t\}$, so we are trying to find $\P(T\in S)$.  Fix $s\in S$ and let $v_s$ be the first vertex along the depth first walk of $s$ whose out-degree is in $A$.  Let $(w_1,\dots, w_d)$ be the elements of $\rspine(v)$ listed by order of appearance on the depth first walk of $s$.  For $1\leq i\leq d$, let $s_i$ be the subtree of $s$ rooted at $w_i$.  Since $\check s= t$, there is a set $J\subseteq \{1,\dots, d\}$ of size $r$ such that $s_i$ contains a vertex with out-degree in $A$ if and only if $i\in J$.  Let $\rho$ be the unique increasing bijection from $\{1,\dots, r\}$ to $J$.  Then, again since $\check s = t$, it follows from the recursive nature of the construction of $\check s$ that $t_i = \check s_{\rho(i)}$.  It follows from these considerations that
\[\begin{split} \P(T\in S) & = \sum_{d=r}^\infty \sum_{s\in S : \#\rspine(v_s)=d}\P(T=s)\\
& = \sum_{d=r}^\infty \P(\rspine(v_T)=d) {d \choose r} \P(N> \tau_{-1})^{d-r}\prod_{j=1}^r\P(\check T' = t_j)
.\end{split}\] 
Since $t_i$ has fewer vertices than $t$, it follows from the induction hypothesis that
\[\P(\check T' = t_i) = \P(N\leq\tau_{-1}) \P(\check T = t_i) = \P(N\leq\tau_{-1})\prod_{v\in t_i} \P(Y = \deg(v)).\]
Thus
\[\begin{split}\P(\check T=t) &= \left[\prod_{j=1}^r\prod_{v\in t_j} \P(Y = \deg(v)) \right] \sum_{d=r}^\infty \P(\rspine(v_T)=d) {d \choose r} \P(N> \tau_{-1})^{d-r} \P(N\leq \tau_{-1})^r\\
& =  \left[\prod_{j=1}^r\prod_{v\in t_j} \P(Y = \deg(v)) \right]\P(Y=r)\\
& = \prod_{v\in t} \P(Y= \deg(v))
,\end{split}\]
which shows that $\check T$ is a Galton-Watson tree with the appropriate offspring distribution.
\end{proof}

Note that $T$ has $n$ vertices with out-degree in $A$ if and only if $\check T$ has $n$ vertices.  This allows us to deduce the following corollary.

\begin{corollary}\label{corollary forest counting}
Let $\F_{n}$ be an ordered forest of $n$ independent Galton-Watson trees with critical or subcritical offspring distribution $\xi$.  Let $(Y_i)_{i\geq 0}$ be distributed like $Y$ from Theorem \ref{theorem check transform} and let $S_i = (Y_1-1)+\cdots+(Y_i-1)$.  Also, let $\F^Y_n$ an ordered forest of $n$ independent Galton-Watson trees with offspring distribution $Y$ and let $Z_n$ be independent with distribution $\textrm{Bin}(n,\P(N\leq \tau_{-1}))$.  Then
\[ \begin{split} \P(\#_A\F_n = k) & = \sum_{j=0}^n \frac{j}{k} {n \choose j}\P(N\leq \tau_{-1})^j\P(N>\tau_{-1})^{n-j} \P(S_k = -j)\\
&= \P(\# \F^Y_{Z_n}=k)
\end{split}.\]
\end{corollary}

\begin{proof}
The idea is to choose a set of $j$ trees to be the trees that contain a vertex with out-degree in $A$ and then apply the Otter-Dwass formula to obtain the formula for the number of vertices in a forest of $j$ independent Galton-Watson trees distributed like $\check T$.  One then sums over the possible sets of $j$ trees. 
\end{proof}

\subsection{The partition at the root}\label{subsec root}
Let $\xi = (\xi_i)_{i\geq 0}$ be a probability distribution with mean $1$ and variance $0<\sigma_1^2<\infty$.  Let $T$ be a Galton-Watson tree with offspring distribution $\xi$ (denote the law of $T$ by $\gw_\xi$).  Let $A\subseteq \Z^+$ and construct $\check T$ as above.  By  Theorem \ref{theorem check transform}, $\check T$ is a Galton-Watson tree.  Let $\zeta$ be its offspring distribution.  Again by Theorem \ref{theorem check transform}, $\zeta$ has mean $1$ and variance $\sigma^2=\gamma^2\sigma_1^2/\xi(A)$, where $\gamma\deq \P(\#_AT\geq 1)$ (see Lemma \ref{lemma Y moments} for computations of the moments of $\zeta$).  Assume that for sufficiently large $n$, $\P(\#_AT=n)>0$.  Let $T^A_n$ be $T$ conditioned to have exactly $n$ vertices with out-degree in $A$ (whenever this conditioning makes sense).

For a $t$ be rooted unordered tree with exactly $n$ vertices with out-degree in $A$, let $\Pi^A(t)$ be the partition of $n$ or $n-1$ (depending on whether or not the degree of the root of $t$ is in $A$) defined by the number of vertices with out-degree in $A$ in the subtrees of $t$ attached to the root.  We also adopt the convention $\Pi^A(\bullet)=\emptyset$, that is, the partition at the root of the tree with one vertex is the emptyset.

\begin{lemma}\label{le root stuff}
(i) Considered as an unordered tree, the law of $T^A_n$ is equal to $\bP^q_n$ where $q_0$ is the law of $T^A_0$ if $0\notin A$ and for $n\geq 1$ such that $T^A_n$ is defined, and $\lambda=(\lambda_1,\dots,\lambda_p)\in \bar \calP^A_n$, we have
\[ q_n(\lambda)=\P(\Pi^A(T^A_n) =\lambda)= \frac{p!}{\prod_{j\geq 1} m_j(\lambda)!} \xi(p) \frac{\prod_{i=1}^p \P(\#_AT= \lambda_i)}{\P(\#_AT = n)} .\]
 \\
(ii) Let $X_1,X_2,\dots$ be i.i.d.\!\! distributed like $\#_AT$ and $\tau_k = X_1+\cdots+X_k$.  We have
\[\P(p(\Pi^A(T^A_n)) =p) = \xi(p)\frac{\P(\tau_p=n-\cf(p\in A))}{\P(\tau_1=n)},\]
and $ \P(\Pi^A(T^A_n) \in \cdot \ | \ \{p(\Pi^A(T^A_n))=p\})$ is the law of a non-increasing rearrangement of $(X_1,\dots,X_p)$ conditionally on $X_1+\cdots+X_p=n-\cf(p\in A)$.
\end{lemma}

\begin{proof}
(i) Letting $c_\emptyset(T^A_n)$ be the root degree of $T^A_n$  and $a_1,\dots,a_p\in \N$ with sum $n-\cf(p\in A)$ we have
\begin{equation} \label{eq root law}  \P(c_\emptyset(T^A_n) = p, \#_A[(T^A_n)_i]=a_i, 1\leq i\leq p) = \xi(p) \frac{\prod_{i=1}^p \P(\#_AT= a_i)}{\P(\#_AT = n)}.\end{equation}
Part (i) now follows by considering the number of sequences $(a_1,\dots, a_p)$ with the same decreasing rearrangement.

(ii) This follows from Equation \eqref{eq root law}.
\end{proof}

To simplify notation, let $q_n$ be the law of $\Pi^A(T^A_n)$ and let $\cf_p=\cf(p\in A)$.  Let $(S_r,r\geq 0)$ be a random walk with step distribution $(\zeta_{i+1}, i\geq -1)$.  Furthermore, since $\gamma=\P(\#_AT\geq 1)$.  By Corollary \ref{corollary forest counting}, we have
\begin{equation}\label{eq A-fund relation}
 q_n(p(\lambda)=p)  = \frac{n\xi(p)}{(n-\cf_p)\gamma \P(S_{n}=-1)}\sum_{j=0}^p {p\choose j}\gamma^j(1-\gamma)^{p-j} j \P(S_{n-\cf_p}=-j)
\end{equation}
where $\hat\xi(p)=p\xi(p)$ is the size-biased distribution of $\xi$.

Define $\bar q_n$ to be the pushfoward of $q_n$ onto $\S^\downarrow$ by the map $\lambda\mapsto \lambda/\sum_i\lambda_i$.

For a sequence $(x_1,x_2,\dots)$ of non-negative numbers such that $\sum_i x_i <\infty$, let $i^*$ be a random variable with
\[\P(i^*=i) = \frac{x_i}{\sum_{j\geq 1} x_j}.\]
The random variable $x^*_1=x_{i^*}$ is called a size-biased pick from $(x_1,x_2,\dots)$.  Given $i^*$, we remove the $i^*$'th entry from $(x_1,x_2,\dots)$ and repeat the process.  This yields a random re-ording $(x^*_1,x^*_2,\dots)$ of $(x_1,x_2,\dots)$ called the size-biased order (if ever no positive terms remain, the rest of the size-biased elements are $0$).  Similarly for a random sequence $(X_1, X_2,\dots)$ we define the size-biased ordering by first conditioning on the value of the sequence.  For any non-negative measure $\mu$ on $\S^\downarrow$, define the size-biased distribution $\mu^*$ of $\mu$ by
\[ \mu^*(f) = \int_{\S^\downarrow} \mu(d\mathbf{s}) \E[f(\mathbf{s}^*)],\]
where $\mathbf{s}^*$ is the size-biased reordering of $\mathbf{s}$.

Define the measure $\nu_2$ on $\S^\downarrow$ by
\[\int_{\S^\downarrow} \nu_2(d\mathbf{s})f(\mathbf{s})= \int_{1/2}^1\sqrt{\frac{2}{\pi s_1^3(1-s_1)^3}} ds_1 f(s_1,1-s_1,0,0,\dots).\]

\begin{theorem}\label{th root limit} With the notation above,
\[\lim_{n\to\infty} n^{1/2}(1-s_1)\bar q_n(d\mathbf{s})= \frac{\sigma_1\sqrt{\xi(A)}}{2}(1-s_1)\nu_2(d\mathbf{s}),\]
where the limit is taken in the sense of weak convergence of finite measures.
\end{theorem}

\begin{proof}
We follow the reductions in Section 5.1 of \cite{HaMi10}.  By Lemma 16 in \cite{HaMi10} (which is a easy variation of Proposition 2.3 in \cite{Bert06}) it is sufficient to show that 
\[\lim_{n\to\infty} n^{1/2}((1-s_1)\bar q_n(d\mathbf{s}))^*= \frac{\sigma_1\sqrt{\xi(A)}}{2}((1-s_1)\nu_2(d\mathbf{s}))^*.\]
Note that for any finite non-negative measure $\mu$ on $\S^\downarrow$ and non-negative continuous function $f:\S_1\to \R$ we have
\[((1-s_1)\mu(d\mathbf{s}))^*(f) = \int_{\S_1} \mu^*(d\mathbf{x})(1-\max\mathbf{x})f(\mathbf{x}).\]
Consequently the theorem follows from the following Proposition.
\end{proof}

\begin{proposition}\label{pr root limit}
Let $f:\S_1\to \R$ be continuous and let $g(\mathbf{x})=(1-\max\mathbf{x})f(\mathbf{x})$.  Then
\[ \sqrt{n}\bar q^*_n(g) \to \frac{\sigma_1\sqrt{\xi(A)}}{\sqrt{2\pi}}\int_0^1\frac{dx}{x^{1/2}(1-x)^{3/2}}g(x,1-x,0,\dots).\]
\end{proposition}

First note that, by linearity, we may assume that $f\geq 0$ and $||f||_\infty \leq 1$.  We begin the proof of this Proposition with several Lemmas regarding the concentration of mass of $\bar q^*_n$.  

\begin{lemma}
\label{le bound 1}
For every $\epsilon>0$, $\sqrt{n}q_n(p(\lambda)>\epsilon\sqrt{n}) \to 0$ as $n\to\infty$.
\end{lemma}

\begin{lemma} \label{le bound 2} For $g$ as in Proposition \ref{pr root limit} we have
\[\Lim{\eta\downarrow 0}\limsup_{n\to\infty} \sqrt{n}\bar q^*_n(|g|\cf_{\{x_1>1-\eta\}}) = 0 \quad \textrm{and}\quad \lim_{n\to\infty} \sqrt{n}\bar q^*_n(\cf_{\{x_1<n^{-7/8}\}}) = 0.\]
\end{lemma}

\begin{lemma} \label{le bound 3}
For every $\eta>0$ we have
\[ \Lim{n\to\infty} \sqrt{n}\bar q^*_n(\cf_{\{x_1+x_2<1-\eta\}})=0.\]
\end{lemma}

\begin{lemma} \label{le integral limit}
There exists a function $\beta_\eta=o(\eta)$ as $\eta \downarrow 0$ such that
\[\begin{split} \Lim{\eta\downarrow 0}\liminf_{n\to\infty} \sqrt{n}\bar q^*_n(g\cf_{\{x_1<1-\eta,x_1+x_2>1-\beta_\eta\}}) & =  \Lim{\eta\downarrow 0}\limsup_{n\to\infty} \sqrt{n}\bar q^*_n(g\cf_{\{x_1<1-\eta,x_1+x_2>1-\beta_\eta\}})  \\
& = \frac{\sigma_1\sqrt{\xi(A)}}{\sqrt{2\pi}}\int_0^1 \frac{g((x,1-x,0,\dots))}{x^{1/2}(1-x)^{3/2}} dx
.\end{split}\]
\end{lemma}

These lemma's are generalizations of the lemmas in Section 5.1 in \cite{HaMi10} to our current setting and the proofs are essentially the same.   We refer the reader to \cite{HaMi10} for the proofs of Lemmas \ref{le bound 1}, \ref{le bound 2}, and \ref{le bound 3}.  We include the proof of Lemma \ref{le integral limit} because it makes clear how the factor of $\sqrt{\xi(A)}$ appears in the scaling limit.

\begin{proof}
Fix $\eta>0$ and suppose that $\eta'\in (0,\eta)$.  Using Lemmas \ref{le bound 1} and \ref{le bound 2} we decompose according to the events $\{p(\lambda)>\epsilon\sqrt{n}\}$ and $\{\mathbf{x}:x_1\leq n^{-7/8}\}$ to get
\begin{multline}\label{eq decomp} \displaystyle \sqrt{n}\bar q^*_n(g\cf_{\{x_1<1-\eta,x_1+x_2>1-\eta'\}})  = o(1)+\sqrt{n}\sum_{1\leq p\leq \epsilon n^{1/2}} q_n(p(\lambda)=p) \\
\displaystyle \times \sum_{n^{1/8}\leq m_1\leq (1-\eta)(n-\cf_p) \atop (1-\eta')(n-\cf_p)\leq m_1+m_2\leq n-\cf_p} \E[g((m_1,m_2,X^*_3,\dots,X^*_p,0,\dots)/(n-\cf_p)) | \tau_p=n-\cf_p,X^*_1=m_1,X^*_2=m_2] \\
\displaystyle \times \frac{pm_1}{n-\cf_p}\frac{\frac{(p-1)m_2}{n-\cf_p}}{1-\frac{m_1}{n-\cf_p}} \P(X_1=m_1)\P(X_2=m_2)\frac{\P(\tau_{p-2}=n-m_1-m_2-\cf_p)}{\P(\tau_p=n-\cf_p)}
.\end{multline}

Observe that, if $1\geq x_1+x_2\geq 1-\eta'$ and $x_1\leq 1-\eta$, then $x_2/(1-x_1)\geq 1-\eta'/\eta$ and $(1-x_1)/x_2\geq 1$.  

Using the local limit theorem and formula \ref{eq A-fund relation} we observe that $n$ and small $\epsilon$, we have
\[ 1-\eta \leq \frac{q_n(p(\lambda)=p)}{\hat\xi(p)} \leq 1+\eta,\]
And similarly for sufficiently large $n$ and small $\epsilon$, we have
\[ \frac{1-\eta}{\sigma\sqrt{2\pi}} \leq (\gamma p)^{-1}n^{3/2} \P(\tau_p=n) \leq \frac{1+\eta}{\sigma\sqrt{2\pi}},\]
for all $1\leq p\leq \epsilon n^{1/2}$.  We note in particular that $\tau_1=X_1=_dX_2$.  Furthermore, for $n^{1/8}\leq m_1\leq (1-\eta)n$ and $m_1+m_2\geq (1-\eta')n$ we have $m_2\geq (\eta-\eta')n$ so that $m_1$ and $m_2$ go to infinity as $n$ does.  Thus, for large $n$ (how large now depends on $\eta'$) we have
\[ \frac{(1-\eta)^2}{2\pi\sigma^2} \leq \gamma^{-2}(m_1m_2)^{3/2} \P(X_1=m_1)\P(X_2=m_2) \leq \frac{(1+\eta)^2}{2\pi\sigma^2}.\]
Now, recall that $f$ is uniformly continuous on $\S_1$.  Furthermore, on the set $\{\mathbf{x}\in \S_1 : x_1+x_2>3/4\}$ we have $\max\mathbf{x} = x_1\vee x_2$ and $\mathbf{x}\mapsto \max\mathbf{x} $ is thus uniformly continuous on this set.  Therefore for $\eta'<(1/4)\wedge \eta^2$ sufficiently small we have 
\[ |g((m_1,m_2,m_3,\dots)/n) - g((m_1,n-m_1,0,\dots)/n)|\leq \eta,\]
for every $(m_1,m_2,\dots)$ with sum $n$ sufficiently large such that $m_1+m_2\geq (1-\eta')n$.  Take $\beta_\eta\deq \eta'$.  

Given the symmetry of the bounds we have just established it is easy to see that the proofs for the $\limsup$ and $\liminf$ will be nearly identical, one using the upper bounds and the other the lower.  We will only write down the proof for the $\liminf$.  For sufficiently large $n$ we have that, up to addition of an $o(1)$ term, $\sqrt{n}\bar q^*_n(g\cf_{\{x_1<1-\eta,x_1+x_2>1-\eta'\}}) $ is bounded below by
\begin{multline*}  \frac{(1-\eta)^3(1-\eta'/\eta)}{(1+\eta)} \sum_{1\leq p\leq \epsilon n^{1/2}}  (p-1)\hat\xi(p) \frac{1}{n-\cf_p}\\
\times  \sum_{n^{1/8}\leq m_1\leq (1-\eta)(n-\cf_p)} (g((m_1,n-m_1-\cf_p,0,\dots)/(n-\cf_p))-\eta)\\
\times  \frac{1}{(m_1/(n-\cf_p))^{1/2}}\frac{1}{(1-m_1/(n-\cf_p))^{3/2}}\frac{1}{\sigma\sqrt{2\pi}}\\
\times  \sum_{ (1-\eta')(n-\cf_p)-m_1\leq m_2\leq n-m_1-\cf_p}\P(\tau_{p-2}=n-m_1-m_2-\cf_p) 
.\end{multline*}
Observe that this last sum is equal to $\sum_{m=0}^{\eta' (n-\cf_p)} \P(\tau_{p-2}=m)$.  By the local limit theorem, this can be made arbitrarily close to $1$ independent of $1\leq p\leq n^{1/2}$.  Using the convergence of Riemann sums (again care must be taken since the integral we get is improper), we have
\begin{multline*}
\liminf_{n\to\infty} \sqrt{n}\bar q^*_n(g\cf_{\{x_1<1-\eta,x_1+x_2>1-\eta'\}}) \\
\geq  \frac{(1-\eta)^3(1-\eta'/\eta)}{1+\eta}  \sum_{p=1}^\infty \gamma (p-1)\hat\xi(p) \int_0^{1-\eta} \frac{dx}{\sigma\sqrt{2\pi}x^{1/2}(1-x)^{3/2}}(g(x,1-x,0,\dots)-\eta)
\end{multline*}
Letting $\eta\downarrow 0$ coupled with observing that $\sum_{p=1}^\infty (p-1)\hat\xi(p)=\sigma_1^2$ and recalling that $\sigma^2=\gamma^2\sigma_1^2/\xi(A)$ completes the proof.
\end{proof}

\begin{proof}[Proof of Proposition \ref{pr root limit}]
Observe that 
\[|\bar q_n^*(g) - \bar q^*_n(g\cf_{\{x_1<1-\eta,x_1+x_2>1-\eta'\}})| \leq \bar q_n^*(|g|\cf_{\{x_1\geq 1-\eta\}}) +  \bar q_n^*(|g|\cf_{\{x_1+x_2\leq 1-\eta'\}}).\]
Fix $\epsilon>0$ and apply Lemmas \ref{le bound 2} and \ref{le integral limit} to find $\eta,\eta'$  such that
\[\sqrt{n}  \bar q_n^*(|g|\cf_{\{x_1\geq 1-\eta\}})  <\frac{\epsilon}{2}\]
and 
\[\left|\sqrt{n} \bar q^*_n(g\cf_{\{x_1<1-\eta,x_1+x_2>1-\beta_\eta\}}) - \frac{\sigma_1\sqrt{\xi(A)}}{\sqrt{2\pi}}\int_0^1\frac{g(x,1-x,0,0,\dots)}{x^{1/2}(1-x)^{3/2}} dx\right| \leq \frac{\epsilon}{2},\]
for large enough $n$.  For this choice of $\eta,\eta'$ and large $n$ we have
\[\left|\sqrt{n}\bar q_n^*(g) - \frac{\sigma_1\sqrt{\xi(A)}}{\sqrt{2\pi}}\int_0^1\frac{g(x,1-x,0,0,\dots)}{x^{1/2}(1-x)^{3/2}} dx\right| \leq \epsilon + \sqrt{n}\bar q_n^*(|g|\cf_{\{x_1+x_2\leq 1-\eta'\}}).\]
By Lemma \ref{le bound 3} the upper bound goes to $\epsilon$ as $n\to\infty$, and the result follows.
\end{proof}

As an immediate corollary of these results, we also identify the unnormalized limit of $\bar q_n$.

\begin{corollary}
$\bar q_n \overset{d}{\rightarrow} \delta_{(1,0,0,\dots)}$.
\end{corollary}

\begin{proof}
Taking $f\equiv 1$ in Proposition \ref{pr root limit} gives $\bar q_n(1-s_1) \to 0$.  Since $L^1$ convergence implies convergence in probability, it follows that for all $0<\eta<1$ we have $\bar q_n(s_1\geq \eta) \to 1$. 
\end{proof}

Note that, as a consequence of Equation \eqref{eq A-fund relation}, we have $q_n(p(\lambda)=p)) \to \hat\xi(p)$.  Thus, while the degree of the root vertex may be large, only one of the trees attached to the root will have noticeable size.

\subsection{Convergence of Galton-Watson trees}  We are now prepared to prove Theorem \ref{th main result}, which, after all of our work above, is a rather straightforward.  The hardest work that remains is to show that 
\[d_{\textrm{GHP}}\left(\frac{1}{\sqrt{n}} T^A_n, \frac{1}{\sqrt{n}} \tilde T^A_n\right) \to 0\]
when $0\notin A$ (see Theorem \ref{th general A} for the definition of $\tilde T^A_n$).  This is accomplished by the next two lemmas.

\begin{lemma}\label{lemma vertices in gw trees}
Let $T^A_n$ be a Galton-Watson tree conditioned to have exactly $n$ vertices with out-degree in $A$.  Then 
\[\frac{\#_{\Z^+}T^A_n}{n} \to \frac{1}{\xi(A)}\]
in probability.
\end{lemma}

\begin{proof}
Let $\mathbf{X}=(X_n)_{n\geq 1}$ be an i.i.d sequence with distribution $\xi$ and let 
\[\tau_{-1}(\mathbf{X}) = \inf\left\{n : \sum_{k=1}^n (X_k-1)=-1\right\}.\]
Further define $N_k(\mathbf{X})$ to be index of the $k$'th occurrence in $\mathbf{X}$ of an element in $A$.  Using the bijection between trees and their depth first queues, we see that $\#_{\Z^+}T^A_n$ is distributed like $\tau_{-1}(\mathbf{X})$ conditioned on $\left( N_n(\mathbf{X}) \leq \tau_{-1}(\mathbf{X}) < N_{n+1}(\mathbf{X})\right)$.  By Corollary \ref{corollary forest counting} and the local limit theorem, we see that 
\[\P(N_n(\mathbf{X}) \leq \tau_{-1}(\mathbf{X})< N_{n+1}(\mathbf{X})) \sim cn^{-3/2}\]
for some $c>0$.  Furthermore, the large deviations concentration of $N_n(\mathbf{X})$ around $1/\xi(A)$ implies that for every $\epsilon>0$ 
\[\P\left(\left|\frac{\tau_{-1}(\mathbf{X})}{n}-\frac{1}{\xi(A)}\right| > \epsilon , \ N_n(\mathbf{X}) \leq \tau_{-1}(\mathbf{X})< N_{n+1}(\mathbf{X})\right)\]
decays exponentially as $n$ goes to infinity and the claim follows immediately.
\end{proof}

\begin{lemma}\label{lemma pruning}
Let $T$ be a Galton-Watson tree with critical, finite variance offspring distribution $\xi$ and let $T^A_n$ be distributed like $T$ conditioned to have exactly $n$ vertices with out-degree in $A$.  Assume that $0\notin A$ and let $v_1,\dots, v_d$ be the roots of the maximal fringe subtrees of $T^A_n$ that contain no vertices with out-degree in $A$, indexed by order of appearance on the depth first walk of $T^A_n$.  Let $\tilde T^A_n$ be the tree obtained by removing the fringe subtrees rooted at $v_1,\dots, v_d$ from $T^A_n$.  Then
\[ d_{\textrm{GHP}}\left(T^A_n, \tilde T^A_n\right) = O(\log n)\]
in probability.
\end{lemma} 

\begin{proof}
Let $T_1,\dots, T_d$ be the fringe subtrees of $T^A_n$ rooted at $v_1,\dots, v_d$ respectively.  It follows from the elementary properties of Galton-Watson trees that, conditionally given $d$, $T_1,\dots, T_d$ are i.i.d such that for $t\in\T^u_{A,0}$
\[\P(T_1 = t) = \frac{1}{\P(\#_AT =0)}\P(T=t)  = \prod_{v\in t} \P(\#_AT =0)^{\deg(v)-1}\xi_{\deg(v)}.\]
The last equality shows that $T_1$ is distributed like a Galton-Watson tree with offspring distribution $\mu$ given by
\[\mu_{i} = \P(\#_AT =0)^{i-1}\xi_{i}\cf(i\notin A).\]
Since $\xi$ has mean $1$ and $0<\P(\#_AT=0)<1$ we see that $\mu$ has mean strictly less than $1$.  It follows that 
\[\P(\textrm{height}(T_1)> k) \leq \beta^k,\]
where $\beta$ is the mean of $\mu$.  In particular, the distribution of the height of $T_1$ decays exponentially.  By Lemma \ref{lemma vertices in gw trees}, there is some $C$ such that $\P(d\leq Cn)\to 1$ as $n\to \infty$.  Observe that
\[\P\left(\max_{1\leq i\leq d} \textrm{height}(T_i) >x , d\leq Cn\right) \leq \P\left(\max_{1\leq i\leq Cn} Z_i >x\right),\]
where $(Z_i,i\geq 1)$ are i.i.d distributed like $\textrm{height}(T_1)$.  Since the distribution of $\textrm{height}(T_1)$ decays exponentially, $\max_{1\leq i\leq Cn} Z_i $ grows logarithmically.  The result now follows from the fact that
\[ d_{\textrm{GHP}}\left(T^A_n, \tilde T^A_n\right) \leq \max_{1\leq i\leq d} \textrm{height}(T_i) \]
almost surely.
\end{proof}

\begin{proof}[Proof of Theorem \ref{th main result}]
Lemma \ref{le root stuff} shows that $T^A_n$ has law $\bP^q_n$ for a particular choice of $(q_n)_{n\geq 1}$.  Theorem \ref{th root limit} and Lemma \ref{lemma pruning} then show that the hypotheses of Theorem \ref{th general A} are satisfied.
\end{proof}

\section{Acknowledgments} I would like to thank my adviser, Jim Pitman, for suggesting the problem that led to this paper and for helpful conversations and feedback along the way.  I would also like to thank the referee for encouragement to remove the condition that $0\in A$, which appeared in earlier versions of this paper.

\bibliographystyle{plain}
\bibliography{Probability_Bibliography}
\end{document}